\documentclass[11pt, draft]{amsart}

\newtheorem{theorem}{Theorem}[section]

\newtheorem{lemma}[theorem]{Lemma}

\newcommand{\ZZ}{{\mathbb Z}}

\newcommand{\TT}{{\mathbb T}}

\newcommand{\DD}{{\mathbb D}}

\newcommand{\cQ}{{\mathcal Q}}

\newcommand{\cR}{{\mathcal R}}

\newcommand{\La}{\Lambda}

\title[  Asymptotics of orthogonal polynomials ]{ Asymptotics of orthogonal polynomials beyond the scope of Szeg\H{o}'s theorem  }
\author{ F. Peherstorfer, A. Volberg, P. Yuditskii}
  \thanks{Partially supported by NSF grant DMS-0200713}


\begin{document}

\begin{abstract}
First we give here a simple proof of a remarkable result of
Videnskii and Shirokov: let $B$ be a Blaschke product with $n$
zeros, then there exists an outer function $\phi, \phi(0)=1$, such
that $\|(B\phi)'\| \leq C n$, where $C$ is an absolute constant.
Then we apply this result to a certain problem of finding the
asymptotic of orthogonal polynomials.
\end{abstract}

\maketitle

\section{Introduction and the main results}
\label{Intr}

Let us introduce the quantity investigated in Videnskii-Shirokov
paper \cite{VS}.

Let $X$ be a Banach space of functions analytic in the open unit
disk $\DD$, let $E$ be a finite collection of points in $\DD$, and
let $B_E$ denote the Blaschke product vanishing precisely at the
points of $E$. Characterization of zero sets of functions from $X$
is a difficult and important problem.  It is closely related to
the estimation of the quantity

$$
\phi(E,X) := \inf\{\|f\|: f \in X, f=B_EG, G(0)=1\}.
$$
The following sequence plays an important part in understanding
zero sets structure of $X$:

$$
r_n(X) :=\sup\{\phi(E,X): card E =n\}.
$$

Very interesting research is related to asymptotics of
$r_n(l^1_a)$, where $l^1_a$ stands for analytic functions with
absolutely summable Taylor coefficients. In the papers of
Sh\H{a}ffer \cite{Sh}, Gluskin, Meyer, Pajor \cite{GMP} and
Queffelec \cite{Q} it was proved that $r_n(l^1_a)\asymp \sqrt{n}$.
Actually Sh\H{a}ffer, motivated by Van der Waerden, studied the
asymptotics of another quantity, which turns out to be equal to
$r_n(l^1_a)$ (see \cite{Sh}, \cite{GMP}, \cite{N}).  In \cite{N} one can find an easy proof of equality 
of  $r_n(l^1_a)$ and this other quantity
$$
k_n := \sup\{|\det T |\|T^{-1}\|: \|T\|\leq 1\},
$$
where the supremum is taken over all invertible matrices of order
$n$ and over {\it all norms}! But if one fixes the norm of $T$ as
acting in $n$-dimensional dimensional Banach space one can
consider the analog of $k_n$, or the following problem of finding
the best estimate
$$
\|T^{-1}\| \leq \Phi_n(\delta),
$$
where we consider all invertible $T$ in $n$-dimensional Banach
space $A$, such that $\|T\|\leq 1$ and the spectral radius
$r(T^{-1})= \delta$. Those estimates are considered in details in
the paper of Nikolski \cite{N}. In \cite{N} the following variant
of the problem is also considered: the requirement $\|T\|\leq 1$
is replaced by
$$
\|f(T)\| \leq C\|f\|,
$$
for every polynomial $f$. In other words operator $T$ is assumed
to satisfy $A$-functional calculus with respect to function space
(algebra) $A$. It is easy to see that for Banach algebras $A=X$
this brings us to estimates of yet another interesting quantity:
$$
cap_X(E) := \inf\{\|f\|_X: f(0) =1, f|E =0\}.
$$
It is immediate that our quantities $\phi(E,X), cap_X(E)$ are
essentially related to each other.

In \cite{N} the sharp estimates of $cap_{B_{p,q}^s}(E)$ is given.
Here $B_{p,q}^s$ is a Besov class (see definition below). Notice
that for non-integer $s$
$$
B_{\infty,\infty}^s = \La^s,
$$
but for integer $s$ only the strict inclusion holds:
\begin{equation}
\label{incl} \La^s\subset B_{\infty,\infty}^s\,.
\end{equation}

One of the results of \cite{N} formulated in the form convenient
for us is

$$
r_n(B_{\infty,\infty}^s)\asymp n^s\,.
$$

By our previous remark this gives $ r_n(\La^s) \asymp n^s $ for
non-integer $s$. The question for integer $s$ seems to be more
subtle. In fact, the proof of this fact for integers given in the
paper of Videnskii and Shirokov (\cite{VS}, Theorem 3.1)  is quite
involved.

We want to give here a very simple proof. It follows the ideas of
\cite{N}, \cite{VS}, but it is quite short. We treat the case
$s=1$.

\begin{theorem}
\label{s1} Let $B$ be a Blaschke product with $n$ zeros in the
disc $\DD$. Then there exists an outer function $\phi$ satisfying
the following properties:
\begin{equation}
\label{prime} \|(B\phi)'\|_{L^{\infty}(\TT)} \leq C\,n\,.
\end{equation}
\begin{equation}
\label{zero} \phi(0) =1\,.
\end{equation}
\begin{equation}
\label{norm} \|\phi\|_{L^{\infty}(\TT)} \leq A\,.
\end{equation}
Here $C,A$ are absolute constants.
\end{theorem}

\subsection{Classes $B_{p,q}^s$,\,$\La^s$}
\label{classes}

Let $W_n$ are Valee-Poussin kernels defined by their Fourier
coefficients as follows. $W_0(z) = 1+z$, and for $n\geq 1$
$\hat{W_n} = 1$ for $k= 2^n$, $\hat{W_n} = 0$ for $k\notin
(2^{n-1}, 2^{n+1})$, and $\hat{W_n}$ is affine on $[2^{n-1},
2^{n}]$ and on $[2^{n}, 2^{n+1}]$.

Now the class $B_{p,q}^s$ consists of analytic in $\DD$ functions
$f(z) = \sum \hat{f}(n) z^n$ such that the sequence $\{
2^{ns}\|f\star W_n\|_{L^p(\TT)}\}$ belongs to $l^q$.

For integer $s=n$ the class $\La^s$consists of all functions
analytic in $\DD$ for which the following seminorm is finite ??

$$
\|f\|_{\La^n}:= \|f^{(n)}\|_{H^{\infty}}\,.
$$
If $s=n+\alpha$, $0<\alpha<1$, then $f\in \La^s$ means that the
previous seminorm is finite and in addition to that one more
seminorm should be finite:

$$
\|f^{(n)}\|_{\La^{\alpha}}:=\sup_{z,\zeta\in \DD,
z\neq\zeta}\frac{|f^{(n)}(z)
-f^{(n)}(\zeta)|}{|z-\zeta|^{\alpha}}\,.
$$

\section{Proof of Theorem \ref{s1}}

\begin{proof}

Let $\widetilde{B}$ be a Blaschke product with the same zeros
$\{z_1,..,z_n\}$ as $B$ but in $(1+\frac{1}{n})\DD$. If
$R=1+\frac{1}{n}$ then

$$
\widetilde{B} = \Pi_{k=1}^n \frac{\frac{z}{R}- \frac{z_k}{R}}{1
-\frac{1}{R^2}\bar{z_k} z}\,.
$$

We can write

$$
R^n \widetilde{B}=: B\phi_0,
$$
having this  as a definition of $\phi_0$.

So

$$
\phi_0(z) = \Pi_{k=1}^n \frac{1- \bar{z_k}z}{1
-\frac{1}{R^2}\bar{z_k} z}\,.
$$
Obviously, if $z\in\TT$ then $|\widetilde{B}| < 1$, and,
therefore,
$$
|\phi_0(z)| \leq (1+\frac{1}{n})^n\leq A\,.
$$
So \eqref{norm} is satisfied. Also $\phi_0(0) = 1$. To check
\eqref{prime} we write the Cauchy formula for
$\widetilde{B}(\zeta_0), \zeta_0 \in \TT$:
$$
\widetilde{B}(\zeta_0) = \frac{1}{2\pi
i}\int_{(1+\frac{1}{n})\TT}\frac{\widetilde{B}(\zeta)}{\zeta
-\zeta_0}d \zeta\,.
$$

Hence

$$
|\widetilde{B}'(\zeta_0)| \leq  \frac{1}{2\pi
}\int_{(1+\frac{1}{n})\TT}\frac{|\widetilde{B}(\zeta)|}{|\zeta
-\zeta_0|^2}d m(\zeta)\,.
$$

And therefore

$$
|\widetilde{B}'(\zeta_0)| \leq
C\,\int_{(1+\frac{1}{n})\TT}\frac{1}{|\zeta -\zeta_0|^2}d
m(\zeta)\leq C\,n\,.
$$

But then
$$
|(B\phi_0)'(\zeta_0)| = (1+\frac{1}{n}))^n
|\widetilde{B}'(\zeta_0)| \leq AC\,n\,.
$$
\end{proof}

\noindent{\bf Remarks}. 1. Replacing $R = 1+\frac{1}{n} $ by $ R =
1+\frac{\varepsilon}{n}$ we can improve the estimate of $\phi_0$:
$\|\phi_0\|_{\infty} \leq 1 + \varepsilon$ by paying the price in
\eqref{prime}, the constant begins to depend on $\varepsilon$:
\begin{equation}
\label{primeepsilon} \|(B\phi)'\|_{L^{\infty}(\TT)} \leq
C_{\varepsilon}\,n\,.
\end{equation}

\noindent 2. We proved also that

$$
\|B\phi_0\|_{\La^s} \leq C\, n^s\,.
$$
For $s\in (0,1)$ this follows immediately from \eqref{prime} and
\eqref{norm}. For $s= 1$ we just presented the proof. For $s>1$
the proof goes exactly along the same lines but we need to
differentiate more times.

\noindent The second remark contains the claim of Theorem 3.1 of
\cite{VS}.

\noindent The first remark plays a crucial part in this paper,
 where we consider the asymptotics of
polynomials orthogonal with respect to measure $\mu$ of the
following type
$$
d\mu := \frac{dm}{|\psi(e^{i\theta})|^2} + \sum_{k=1}^{\infty}
\mu_k \delta_{z_k}\,,
$$
where $z_k$ form the sequence of Blaschke.

In this situation one should distinguish between the sequence of
Blaschke inside $\DD$ and outside $\DD$. The second case is much
more interesting as expected asymptotics of $L^2(\mu)$-orthogonal
polynomials is different from the classical asymptotics of
Szeg\H{o}, where all $\mu_k =0$.

The orthogonalization with respect to measures described above is
an attempt to generalize the results of the paper of Peherstorfer
and Yuditskii \cite{PY}, where the "point spectrum" was assumed to
be lying on the real line. This is  quite an interesting paper
which extends the strong asymptotic results for orthogonal
polynomials that satisfy a Szeg\H{o} condition on a real interval
to the case where the measure has in addition a denumerable set of
mass points outside the interval. The authors assume that the mass
points outside of the interval accumulate only at the endpoints of
the interval and that a Blaschke condition is satisfied. Under
these conditions, strong asymptotic results are given for the
orthonormal polynomials and for their leading coefficients. The
asymptotic behavior is stated in terms of the Szeg\H{o} function
for the absolutely continuous part and of a Blaschke product
related to the extra mass points.

\section{Orthogonal polynomials and smoothness of $B\phi$}
\label{orth}

Let $\mu$ be as above with
$$
|z_k|>1, \,\, \sum_k (|z_k| -|) < \infty,
$$
and let $\{P_n\}$ denote the sequence of analytic polynomials
orthonormal with respect to $\mu$.

We are writing
$$
P_n(z) = \tau_n z^n +...+ a_{0n}
$$
and we want to find the limit of $\tau_n$ (we will prove that it
exists in many cases).

We also consider rational functions of the form

\begin{equation}
\label{Rnm}
R_{n,m}(z) = \eta_{n, m} z^n + ...+ a_{0n} +...+ a_{m,n}z^{m}\,,
\end{equation}
where $m = -n$ or $-(n-1)$. We want them to be mutually orthogonal in
$L^2(\mu)$ and of norm $1$ in this space. We wish to discuss
the asymptotics of the ``leading" coefficient $\eta _n := \eta_{n, -( n-1)}$.

Here are our two main results. We have the convention that
$\tau_n, \eta_n = \eta_{n, -(n-1)}$ are all positive, and that $B(0), \psi(0)$ are
also positive. The first theorem deals with the perturbation of so-called CMV matrices.

\begin{theorem}
\label{eta} For any $\{z_k\}$ satisfying the Blaschke condition as
above, for any summable $\{\mu_k\}$, and for $\psi $ bounded away
from $0$  one has
$$
\lim_{n\rightarrow\infty} \eta_n = B(0)\psi(0)\,.
$$
\end{theorem}

The asymptotics of $\tau_n$ seems to be more subtle. For certain
geometric configuration of $\{z_k\}$   no extra assumptions on
$\{\mu_k\}$ is needed.  For example, if $\{z_k\}$ converge in a
Stolz star to only a closed set $E\subset \TT$ satisfying $m(E)=0$
and $\sum l_n\log\frac{1}{l_n} <\infty$ (here $l_n$ stand for
lengths of complimentary intervals of $E$; such sets are called
Carleson subsets of $\TT$) we do not need any extra assumption on
$\{\mu_k\}$. A particular case of such a geometry is the case of
real $\{z_k\}$ considered in \cite{PY}. But if $\{z_k\}$ is an
arbitrary Blaschke sequence outside of $\DD$, we still need one
assumption on $\{\mu_k\}$ saying that the series $\sum_k\mu_k$
converges with a definite speed. Here is this assumption:

\begin{equation}
\label{log} \sum_{1< |z_k| < 1 + \frac{1}{n}} \mu_k \leq C_A
\frac{1}{(\log n)^A}, \,\, \forall \, A\,.
\end{equation}

\begin{theorem}
\label{tau} For any $\{z_k\}$ satisfying the Blaschke condition as
above, for any summable $\{\mu_k\}$ satisfying the extra condition
\eqref{log}, and for $\psi $ bounded away from $0$ and from
$\infty$ one has
$$
\lim_{n\rightarrow\infty} \eta_n = B(0)\psi(0)\,.
$$
\end{theorem}

\section{The proof of Theorem \ref{eta}}
\label{proofeta}

\begin{lemma}
\label{extremaleta}
Let $R_n$ be a rational function solving the following extremal problem: find
$$
\sup\{|\eta|: \cR(z) = \eta z^n +...+ a z^{-(n-1)}, \|R\|_{L^2(\mu)}\leq 1\}.
$$
Then $R_n := R_{n, -(n-1)}$ from \eqref{Rnm} solves this extremal problem.
\end{lemma}

\begin{proof}

Consider another extremal problem: $\inf\{\|Q\|_{L^2(\mu)}:  \cQ(z) =  z^n +...+ a z^{-(n-1)}\}$. These are
$L^2(\mu)$ orthogonal functions of course. Our extremal functions $\cR$ are just  $\cQ/\|\cQ\|_{L^2(\mu)}$, so $\inf \|Q\|_{L^2(\mu)}$ means $\sup |\tau|$.
\end{proof}

\subsection{Estimate from below of $\eta_n$}
\label{beloweta}

Let $\zeta_k = 1/\bar{z_k}$, $B$ be a Blaschke product
with zeros $\{\zeta_k\}_{k=1}^{\infty}$, $B_n$ be a partial Blaschke product built by
$$
\zeta_k:  |\zeta_k| < 1 - \frac{1}{k_n},
$$
where $k_n := [\varepsilon_n n]$, $\varepsilon\rightarrow 0$.

We choose  $\varepsilon_n$ will be chosen later.

Without the loss of generality we think that
$$
\sum_k (1 -|\zeta_k|) < 1.
$$
The number of zeros of $B_n$ is at most $k_n$.

Choose $k_n := [A_n\varepsilon_n n]$, where $A_n$ grows to infinity very slow, namely, 
$A_n\varepsilon\rightarrow 0$.
Put
$$
R:= 1 + \frac{1}{l_n}\,.
$$
Let $\widetilde{B}$ be a Blaschke product with the same zeros as $B_n$, but a Blaschke product in $R\DD$. Then we introduce $\phi$ by formula (as in Theorem \ref{s1})
\begin{equation}
\label{phidef}
R^{k_n} \widetilde{B}=: B\phi\,,
\end{equation}
or 
$$
\phi  := \Pi_{k=1}^{k_n} \frac{1- \bar{\zeta_k} z}{1 -\frac{1}{R^2} \bar{\zeta_k} z}\,.
$$

 In Theorem \ref{s1}  we proved that 
 
 $$
 \|B\phi\|_{\La^1} \leq C l_n\,.
 $$
 
 Therefore,
 
 $$
 \|B\phi\|_{B_{\infty,\infty}^1} \leq C l_n\,.
 $$
 
 Let  $2^k$ be the largest such number smaller than $n$.  Consider $W_{k}, W_{k+1},...$
 (we introduced Vall\'ee-Poussin kernels in Section
 \ref{Intr}). We just noticed that
 
 \begin{equation}
 \label{VP1}
\|W_j\star B\phi\|_{\infty} \leq C\, l_n 2^{-j}, \,\, j\geq k\,.
\end{equation}
So \eqref{VP1} gives
 \begin{equation}
 \label{VP2}
\|\sum_{j\geq k}W_j\star B\phi\|_{\infty} \leq C\, l_n/n
\end{equation}
And, therefore, if $V_k$ denote the modified Vall\'ee-Poussin kernel 
$$
\hat{V}_k (j) =1,\, j\in[-2^{k-1}, 2^{k-1}],\,\,\hat{V}_k (j) =1,\, |j|\geq 2^k,
$$
and $\hat{V}_k (j) $ is affine otherwise, then
from \eqref{VP2} we conclude
\begin{equation}
 \label{VP3}
\|V_k\star (B\phi) - B\phi\|_{\infty} \leq C\, l_n/n
\end{equation}

Let us consider the following modified Vall\'ee-Poussin kernel 
$$
\hat{VP}_n(j) =1,\, j\in[-n, n],\,\,\hat{VP}_n (j) =1,\, |j|\geq n,
$$
and $\hat{VP}_n (j) $ is affine otherwise.  It is easy to see that we proved (using $V_k\star VP_n = V_k$)

\begin{equation}
 \label{VP4}
\|VP_n\star (B\phi) - B\phi\|_{\infty} \leq 4C\, l_n/n \leq C'\delta_n\,,
\end{equation}
where $\delta_n \rightarrow 0$.

Notice the estimate from above for $\phi$ (keeping in mind $l_n\approx A_n  k_n$):

\begin{equation}
\label{phiabove}
|(B_n\phi)(z)| = (1 +\frac{1}{l_n})^{k_n} |\widetilde{B}(\zeta)| \leq e^{C/A_n}=:1+ \delta_n^{'},\,\delta_n^{'}\rightarrow 0\,.
\end{equation}

Consider  $G_n = VP_n\star (B\phi), r_n = G_n/z^n$. Clearly (using $\phi(0)=1$)
\begin{equation}
 \label{starcoef}
|G_n(0) - B(0)\phi(0)|  = |G_n(0) - B(0)|\leq \delta_n\,.
\end{equation}

The rational function $r_n$ is of the type we want, but only after the  application of the symmetry $n\rightarrow -n$  on $\ZZ$.

Let us estimate the norm

$$
\|r_n\|^2:= \int_{\TT} |r_n|^2\,dm + \sum_{k=1}^{\infty} \mu_k |r_n(\zeta_k)|^2\,.
$$

From \eqref{VP4} and \eqref{phiabove} we have $|r_n| \leq 1 + C(\delta_n+ \delta_n^{'})$ on $\TT$. So $\int_{\TT} |r_n|^2\,dm \leq 1 + C\delta_n^{''},\,\delta_n^{''}\rightarrow 0$.
Now we split the sum to two: $\alpha_1:=\sum_{k=1}^{k_n} \mu_k |r_n(\zeta_k|^2$, $\alpha_2:=\sum_{k>k_n} \mu_k |r_n(\zeta_k)|^2$. 

To estimate $\alpha_1$ notice that \eqref{VP4} says that $G_n -B\phi$ is small on $\TT$, and the construction of $G_n$ says that $(G_n -B\phi)(z)$  has zero of multiplicity $n$ at the origin. This is analytic function in $\DD$, so by classical Schwartz lemma

\begin{equation}
 \label{VP5}
|G_n(z) - B\phi(z)|\leq C\delta_n\,|z|^n,
\end{equation}
where $\delta_n \rightarrow 0$.

In particular,
\begin{equation}
 \label{VP6}
|(VP_n \star B\phi)(\zeta_k) - B\phi(\zeta_k)|\leq C\delta_n\,|\zeta_k|^n\,.
\end{equation}

But $B(\zeta_k) =0,\, k=1,..., k_n$, so we rewrite

\begin{equation}
 \label{VP7}
|G_n(\zeta_k)|= |(VP_n \star B\phi)(\zeta_k) |\leq C\delta_n\,|\zeta_k|^n,\, k=1,..., k_n\,.
\end{equation}

And finally, using $r_n(z) = G_n(z)/z^n$, we get
\begin{equation}
 \label{VP7}
|r_n(\zeta_k) |\leq C\delta_n,\,\, k=1,..., k_n\,.
\end{equation}

Therefore,
\begin{equation}
 \label{alpha1}
|\alpha_1| \leq C\delta_n \sum_k \mu_k \leq C\delta_n\,.
\end{equation}

To estimate $\alpha_2$ we start from \eqref{VP6} to write
\begin{equation}
 \label{VP9}
|r_n(\zeta_k) |\leq C\delta_n + |(B\phi)(\zeta_k)|/|\zeta_k|^n,\, k>k_n.
\end{equation}
But if $k>k_n$ then $1-|\zeta_k| < \frac{C}{\varepsilon_n n}$, so
$$
\frac{1}{|\zeta_k|^n} \leq e^{C/\varepsilon_n}\,.
$$
In particular,
$$
|\alpha_2|\leq \sum_{k>k_n} \mu_k|r_n(\zeta_k) |^2 \leq C\delta_n +  e^{C/\varepsilon_n}\ \sum_{k>k_n} \mu_k\,.
$$

We can always choose $\varepsilon_n$ to decrease so slowly that $\sum_{k>k_n} \mu_k$ kills the growth of $e^{C/\varepsilon_n}$ and the last expression
$$
\gamma_n :=e^{C/\varepsilon_n}\sum_{k>k_n} \mu_k
$$
tends to zero.

Summarizing all that
\begin{equation}
 \label{etaI}
\|r_n\| \leq 1 + \text{small}
\end{equation}
and
\begin{equation}
 \label{etaII}
 r_n(z) = \rho z^{-n} +...+ a z^{n-1}, |\rho| \geq B(0)  - \text{small}\,.
\end{equation}

Let us inverse the variable. Consider $R(z) := r^*_n(z) := \overline{r(1/\bar{z})}= \bar{\rho} z^{n} +...+ \bar{a} z^{-(n-1)}$. 
Then we just got $|\rho| \geq B(0) \phi(0) - \text{small}$ and
\begin{equation}
 \label{etaI}
\|R\| _{L^2(\mu)}\leq 1 + \text{small}
\end{equation}

We finished the estimate

\begin{equation}
\label{etafrombelow}
\liminf_{n\rightarrow \infty} |\eta_n| \geq B(0)\,.
\end{equation}

This proves the estimate from below for $\eta_n$ for the case $\psi =1$. For non-trivial $\psi$ bounded away from zero, we approximate it from below by smooth $\psi'$ and repeat the above argument
for $B\phi/\psi$ instead of $B\phi$.

\section{Estimate from below of $\tau_n$}
\label{belowtau}

We repeat   the beginning of Section \ref{beloweta}.
Let $\zeta_k = 1/\bar{z_k}$, $B$ be a Blaschke product
with zeros $\{\zeta_k\}_{k=1}^{\infty}$, $B_n$ be a partial Blaschke product built by
$$
\zeta_k:  |\zeta_k| < 1 - \frac{1}{k_n},
$$
where $k_n := [\varepsilon_n n]$, $\varepsilon\rightarrow 0$.

We choose  $\varepsilon_n$ later.

Without the loss of generality we think that
$$
\sum_k (1 -|\zeta_k|) < 1.
$$
The number of zeros of $B_n$ is at most $k_n$.

Choose $k_n := [A_n\varepsilon_n n]$, where $A_n$ grows to infinity very slow, namely, 
$A_n\varepsilon\rightarrow 0$.
Put
$$
R:= 1 + \frac{1}{l_n}\,.
$$
Let $\widetilde{B}$ be a Blaschke product with the same zeros as $B_n$, but a Blaschke product in $R\DD$. Then we introduce $\phi$ by formula (as in Theorem \ref{s1})
\begin{equation}
\label{phidef}
R^{k_n} \widetilde{B}=: B\phi\,,
\end{equation}
or 
$$
\phi := \Pi_{k=1}^{k_n} \frac{1- \bar{\zeta_k} z}{1 -\frac{1}{R^2} \bar{\zeta_k} z}\,.
$$

We use all the notations of the previous section.
Let  $T_n$ denote the Taylor polynomial of degree $n$ of $B_n\phi$.
Our goal is very simple, we want to repeat all the reasoning of the previous section, but instead of estimating $|G_n(z) - B\phi(z)| = |(VP_n \star B\phi)(z) - B\phi(z)|$ we wish to estimate $|T_n(z) - B\phi(z)|=
|(D_n \star B\phi)(z) - B\phi(z)|$. Here $D_n$ stands for Dirichlet kernel. In other words we just want to replace $VP_n$ kernel by $D_n$ kernel.

Let  $\zeta_0\in \TT$. One can write the Cauchy formula for $B\phi(\zeta_0)$ over $R\TT, R = 1 +\frac{1}{l_n}$:
$$
(B_n\phi - T_n)(\zeta_0) =\frac{1}{2\pi i} \int_{(1 +\frac{1}{l_n})\TT}\frac{\zeta_0^n (B_n\phi)(\zeta)}{\zeta - \zeta_0) \zeta^n}d\zeta\,.
$$

We saw in \eqref{phiabove} that

\begin{equation}
\label{phiabove1}
|B_n\phi(\zeta)| = (1 +\frac{1}{l_n})^{k_n} |\widetilde{B}(z)| \leq e^{C/A_n}=: 1+ \kappa_n\,, \kappa_n\rightarrow 0\,.
\end{equation}
In particular, putting this in our Cauchy formula we get

$$
\|B_n\phi - T_n\|_{\infty} \leq \frac{C}{(1 +\frac{1}{l_n})^n} \int_{(1 +\frac{1}{l_n})\TT}\frac{|d\zeta|}{|\zeta- \zeta_0|}\,,
$$
or
$$
\|B_n\phi - T_n\|_{\infty} \leq C\, \log l_n  \,e^{-\frac{1}{A_n\varepsilon_n}}\leq  C'\,\log n\,e^{-\frac{1}{A_n\varepsilon_n}}\,.
$$

We denote $\delta_n := \log n\,e^{-\frac{1}{A_n\varepsilon_n}}$. Unlike $\delta_n$ from the previous section, this one should not go to zero.
But we require that it does go to zero, thus putting  some restriction on $\varepsilon_n, A_n$.  So from now on we {\it assume}
\begin{equation}
\label{assump}
\delta_n:= \log n\,e^{-\frac{1}{A_n\varepsilon_n}\rightarrow 0},
\end{equation}
and
\begin{equation}
\label{smallT}
\|B_n\phi - T_n\|_{\infty} \leq C\delta_n \rightarrow 0\,.
\end{equation}
Also  introduce
$$
 p_n(z) := T_n(z)/z^n,\,P_n(z):=\overline{p_n(\frac{1}{\bar{z}})} = z^n\overline{T(\frac{1}{\bar{z}})}\,.
$$

Then, let us consider the norm
$$
\|p_n\|:= \int_{\TT} |p_n|^2 dm + \sum_k \mu_k |p_n(\zeta_k|^2\,.
$$

Of course we want to prove
\begin{equation}
\label{tau1}
\|p_n\| \leq 1 + \text{small}
\end{equation} 
and
\begin{equation}
\label{tau2}
|T_n(0)| \geq B(0)  -  \text{small}
\end{equation} 
Relation \eqref{tau2} follows from \eqref{smallT} and $\phi(0) =1$.

To prove \eqref{tau1} we notice that \eqref{phiabove1} implies
$$
\int_{\TT}|p_n|^2 dm \leq 1 + C(\kappa_n + \delta_n) \leq 1 + \text{small}\,.
$$

We split the sum $\sum_k \mu_k |p_n(\zeta_k|^2$ into two: $\beta_1 := \sum_{k=1}^{k_n} \mu_k |p_n(\zeta_k|^2$,  $\beta_2 := \sum_{k>k_n} \mu_k |p_n(\zeta_k|^2$.

To estimate $\beta_1$ we notice \eqref{smallT} and the fact that $B\phi - T_n$ is an analytic function in $\DD$ with zero of multiplicity $n$ at the origin.
 So by Schwartz lemma

\begin{equation}
 \label{T5}
|T_n(z) - B\phi(z)|\leq C\delta_n\,|z|^n,
\end{equation}
where $\delta_n \rightarrow 0$.

In particular,
\begin{equation}
 \label{T6}
|T_n(\zeta_k) - B\phi(\zeta_k)|\leq C\delta_n\,|\zeta_k|^n\,.
\end{equation}

But $B(\zeta_k) =0,\, k=1,..., k_n$, so we rewrite

\begin{equation}
 \label{T7}
|T_n(\zeta_k)| \leq C\delta_n\,|\zeta_k|^n,\, k=1,..., k_n\,.
\end{equation}

And finally, using $p_n(z) = T_n(z)/z^n$, we get
\begin{equation}
 \label{T7}
|p_n(\zeta_k) |\leq C\delta_n,\,\, k=1,..., k_n\,.
\end{equation}

Therefore,
\begin{equation}
 \label{beta1}
|\alpha_1| \leq C\delta_n \sum_k \mu_k \leq C\delta_n\,.
\end{equation}

To estimate $\beta_2$ we start from \eqref{T6} to write
\begin{equation}
 \label{T9}
|p_n(\zeta_k) |\leq C\delta_n + |(B\phi)(\zeta_k)|/|\zeta_k|^n,\, k>k_n.
\end{equation}
But if $k>k_n$ then $1-|\zeta_k| < \frac{C}{\varepsilon_n n}$, so
$$
\frac{1}{|\zeta_k|^n} \leq e^{C/\varepsilon_n}\,.
$$
In particular,
$$
|\beta_2|\leq \sum_{k>k_n} \mu_k|p_n(\zeta_k) |^2 \leq C\delta_n +  e^{C/\varepsilon_n} \sum_{k>k_n} \mu_k\,.
$$

We cannot always choose $\varepsilon_n$ to decrease so slowly that $\sum_{k>k_n} \mu_k$ kills the growth of $e^{C/\varepsilon_n}$ and the last expression
\begin{equation}
\label{assump2}
\gamma_n :=e^{C/\varepsilon_n}\sum_{k>k_n} \mu_k
\end{equation}
tends to zero. The reason is in the fact that we have a restriction \eqref{assump} that makes $\varepsilon_n $ to tend to zero with a {\it certain} speed.
But combining \eqref{assump}, \eqref{assump2}, we see that they are reconcilable if $\mu_k$ tend to zero sufficiently rapidly.

It is easy to deduce that
\eqref{log} allows us to have both assumptions \eqref{assump}, \eqref{assump2}. Therefore, we get \eqref{tau1} and \eqref{tau2}. Inversing the variable, that is passing from $p_n(z)$ to $P_n(z)$ proves 
$$
\liminf_{n\rightarrow \infty} |\tau_n| \geq B(0)\,.
$$

If $\psi$ is bounded away from zero we approximate it from below by a smooth function, and thinking now that $1/\psi$ is smooth, we can repeat our approximation arguments with $B\phi/\psi$ instead of $B\phi$.

\section{Estimate from above of $\tau_n$ and $\eta_n$}
\label{above}

This is a simple estimate. For $\eta_n$ it is still slightly more difficult, so we show it in this case.

We want to prove now
\begin{equation}
\label{etaabove}
\limsup_{n\rightarrow \infty} |\tau_n| \leq B(0)\psi(0)\,.
\end{equation}
Let $\psi_*(z) := \overline{\psi(1/\bar{z})}, B_*(z) := \overline{B(1/\bar{z})}$.
Let $B^k_*$ be a finite Blaschke product that is subproduct of $B_*$.

Consider
$$
\int_{\TT} \frac{R_n(z)}{\psi_*(z) }\overline{z^n B^k_*(z)}\, d m(z)\,,
$$ 
which absolute value evidently does not exceed one. Let us calculate the integral according to Cauchy theorem in the exterior of $\DD$. We can do that as $R_n/z^{n+1}$ is analytic in the exterior of the unit disc and equals zero at infinity.

$$
\int_{\TT} \frac{R_n(z)}{\psi_*(z) z^{n +1}B^k_*(z)}\, d z= \frac{\eta_n}{B(0)\psi(0)} - \sum\frac{R_n(z_i)}{(B^k_*)'(z_i)\psi_*(z_i) z_i^{n +1}}\,.
$$
The last sum we have the estimate (since $R_n$ is normalized in $L^2(\mu)$):

$$
\Bigl|\sum\frac{R_n(z_i)}{(B^k_*)'(z_i)\psi_*(z_i) z_i^{n +1}}\Bigr|^2 \leq \sum \frac{1}{|(B^k_*)'(z_i)\psi_*(z_i) z_i^{n +1}|^2 \mu_i}\,.
$$
The number of terms in the sum is finite, so the sum goes to zero when $n$  tends to infinity.

Thus,
$$
\limsup_{n\rightarrow \infty} |\tau_n| \leq B^k_*(0)\psi(0)\,.
$$
Since $B^k_*$ is arbitrary we proved \eqref{etaabove}.

Of course we the reader can see that the proof for $\tau_n$ is exactly the same.

Our both theorems are completely proved.

\end{document}